\documentclass[11pt,a4paper]{amsart}
\usepackage[all]{xy}
\usepackage{amscd}

\newcommand{\ie}{{\em i.e.}\ }
\newcommand{\cf}{{\em cf.}\ }
\newcommand{\eg}{{\em e.g.}\ }
\newcommand{\ko}{\: , \;}
\newcommand{\ul}[1]{\underline{#1}}

\renewcommand{\tilde}[1]{\widetilde{#1}}

\newtheorem*{theorem}{Theorem}
\newtheorem*{lemma}{Lemma}
\newtheorem*{proposition}{Proposition}
\newtheorem*{corollary}{Corollary}

\newcommand{\opname}[1]{\operatorname{\mathsf{#1}}}

\renewcommand{\mod}{\opname{mod}\nolimits}

\newcommand{\Mod}{\opname{Mod}\nolimits}

\newcommand{\ind}{\opname{ind}\nolimits}
\newcommand{\per}{\opname{per}\nolimits}
\newcommand{\add}{\opname{add}\nolimits}

\renewcommand{\ker}{\opname{ker}\nolimits}

\newcommand{\Z}{\mathbb{Z}}

\newcommand{\iso}{\stackrel{_\sim}{\rightarrow}}

\newcommand{\Hom}{\opname{Hom}}
\newcommand{\End}{\opname{End}}
\newcommand{\RHom}{\opname{RHom}}
\newcommand{\HOM}{\opname{Hom^\bullet}}
\newcommand{\Ext}{\opname{Ext}}

\newcommand{\ca}{{\mathcal A}}

\newcommand{\cc}{{\mathcal C}}
\newcommand{\cd}{{\mathcal D}}
\newcommand{\ce}{{\mathcal E}}

\newcommand{\ch}{{\mathcal H}}

\newcommand{\cm}{{\mathcal M}}

\newcommand{\cp}{{\mathcal P}}
\newcommand{\cR}{{\mathcal R}}

\newcommand{\ct}{{\mathcal T}}

\newcommand{\stab}{\opname{stab}}
\newcommand{\cmac}{\opname{CM}}

\let\cal\mathcal
\def\Ascr{{\cal A}}

\def\Cscr{{\cal C}}

\def\Tscr{{\cal T}}

\def\Id{\operatorname{id}}

\def\Lotimes{\overset{L}{\otimes}}

\def\Mod{\operatorname{Mod}}
\def\mod{\operatorname{mod}}

\def\Ext{\operatorname {Ext}}
\def\Hom{\operatorname {Hom}}

\def\End{\operatorname {End}}
\def\RHom{\operatorname {RHom}}

\def\ker{\operatorname {ker}}

\def\End{\operatorname {End}}

\def\add{\operatorname {add}}

\def\r{\rightarrow}

\newtheorem{lemmas}{Lemma}[subsection]

\begin{document}

\title[Acyclic Calabi-Yau categories]{Acyclic Calabi-Yau
categories\\}

\author{Bernhard Keller}
\address{UFR de Math\'ematiques, Institut de Math\'ematiques de Jussieu,
   UMR 7586 du CNRS, Case 7012,
   Universit\'e Paris 7, 2 place Jussieu, 75251 Paris Cedex 05,
   France}
\email{keller@math.jussieu.fr}
\author[Idun Reiten]{Idun Reiten\\
\mbox{ } \\
with an appendix by Michel Van den Bergh}
\address{Institutt for matematiske fag, Norges
Teknisk-naturvitenskapelige universitet, N-7491, Trondheim, Norway}
\thanks{I.~R. supported by a grant from the Norwegian Research Council}
\email{idunr@math.ntnu.no}



\begin{abstract} We prove a structure theorem for triangulated
  Calabi-Yau categories: An algebraic $2$-Calabi-Yau triangulated
  category over an algebraically closed field is a cluster category iff
  it contains a cluster tilting subcategory whose quiver has no
  oriented cycles. We prove a similar characterization for higher
  cluster categories. As an application to commutative algebra, we
  show that the stable category of maximal Cohen-Macaulay modules over
  a certain isolated singularity of dimension three is a cluster
  category. This entails the classification of the rigid Cohen-Macaulay
  modules first obtained by Iyama-Yoshino.
  As an application to the combinatorics of quiver mutation,
  we prove the non-acyclicity of the quivers of endomorphism algebras
  of cluster-tilting objects in the stable categories of
  representation-infinite preprojective algebras. No direct
  combinatorial proof is known as yet. In the appendix,
  Michel Van den Bergh gives an alternative proof of the main theorem
  by appealing to the universal property of the triangulated orbit
  category.
\end{abstract}

\maketitle


\section{Introduction}

Cluster algebras were introduced and studied by Fomin-Zelevinsky and
Berenstein-Fomin-Zelevinsky in a series of articles
\cite{FominZelevinsky02} \cite{FominZelevinsky03}
\cite{BerensteinFominZelevinsky05} \cite{FominZelevinsky07}. It was
the discovery of Marsh-Reineke-Zelevinsky
\cite{MarshReinekeZelevinsky03} that they are closely connected to
quiver representations.  This link is similar to
the one between quantum groups and quiver representations discovered
by Ringel \cite{Ringel90} and investigated by Kashiwara, Lusztig,
Nakajima and many
others.  The link between cluster algebras and quiver representations
becomes especially beautiful if, instead of categories of quiver
representations, one considers certain triangulated categories deduced
from them: the so-called {\em cluster categories}.
These were introduced in
\cite{BuanMarshReinekeReitenTodorov04} and, for Dynkin quivers
of type $A_n$, in \cite{CalderoChapotonSchiffler04}.
If $k$ is a field and $Q$
a quiver without oriented cycles, the associated cluster category $\cc_Q$ is
the `largest' $2$-Calabi-Yau category under the derived category of
representations of $Q$ over $k$. It was shown
\cite{BuanMarshReiten06}
\cite{BuanMarshReitenTodorov05}
\cite{CalderoChapoton06}
\cite{CalderoKeller05a}
\cite{CalderoKeller05b}
that this category
fully determines the combinatorics of the cluster algebra
associated with $Q$ and carries considerably more information.
This was used to prove significant new results on cluster algebras,
\cf \eg
\cite{BuanReiten06}
\cite{CalderoReineke06}.
We refer to
\cite{FominZelevinsky03a} \cite{Zelevinsky04} for more background on
cluster algebras and to
\cite{AssemBruestleSchifflerTodorov05}
\cite{AssemBruestleSchiffler06}
\cite{BuanMarshReiten04b}
\cite{BuanReiten05b}
\cite{BuanReitenSeven06}
\cite{BaurMarsh06}
\cite{GeissLeclercSchroeer06}
\cite{GeissLeclercSchroeer05c}
\cite{Iyama05a}
\cite{IyamaReiten06}
\cite{Ringel07}
\cite{Thomas06}
\cite{Zhu06}
for recent developments in the study of their links with
representations of quivers and finite-dimensional algebras.

The question arises as to whether, for a given quiver $Q$ without
oriented cycles, the cluster category is the `unique model' of the
associated cluster algebra. In other words, if we view the
cluster algebra as a combinatorial invariant associated with
the cluster category, is the category determined by this invariant ?

In this paper, we show that surprisingly, this question has a positive
answer. Namely, we prove that if $k$ is an algebraically closed
field and $\cc$ an
algebraic $2$-Calabi-Yau category containing a cluster tilting
object $T$ whose endomorphism algebra has a quiver $Q$ without
oriented cycles, then $\cc$ is triangle equivalent to the cluster
category $\cc_Q$. Notice that this result is `of Morita type',
but much stronger than typical Morita theorems,
since we only need to know the {\em quiver} of the endomorphism algebra,
not the algebra itself.

We give several applications: First, we show on an example that
cluster categories naturally appear as stable categories of
Cohen-Macaulay modules over certain singularities. This yields an
alternative proof of Iyama-Yoshino's \cite{IyamaYoshino06}
classification of rigid Cohen-Macaulay modules
over a certain isolated singularity.
More examples may be obtained from
\cite{IyamaReiten06} and \cite{IyamaYoshino06}, \cf
\cite{BurbanIyamaKellerReiten07}.

Secondly, we show that the quivers associated in
\cite{GeissLeclercSchroeer05b} with representation-infinite
finite-dimensional preprojective algebras are not mutation-equivalent
to quivers without oriented cycles. This last result was obtained
independently by C.~Geiss \cite{Geiss06}. It has been used in
\cite{DerksenWeymanZelevinsky07}, Example~8.7, to show that the class
of rigid quivers with potential is strictly greater than the class of
quivers with potential mutation-equivalent to acyclic ones.  An
application to the realization of cluster categories as stable
categories of Frobenius categories with finite-dimensional morphism
spaces is given in \cite{BuanIyamaReitenScott07} and
\cite{GeissLeclercSchroeer06a}. In \cite{HolmJoergensen06} and
\cite{HolmJoergensen06a}, the authors use our results to determine
which stable categories of representation-finite selfinjective
algebras of type $A$ and $D$ are higher cluster categories.  More
generally, in \cite{Amiot07}, the author obtains
a classification of `most' triangulated categories with finitely many
indecomposables by methods similar to ours.

The main difficulty in the proof is the construction of a triangle
functor between the cluster category and the given Calabi-Yau
category.  Our construction is based on the description
\cite{Keller05} of the cluster category as a stable derived category of a
certain differential graded category. This approach leads to
interesting connections between Calabi-Yau categories of dimensions
$2$ and $3$, which have been further investigated in
\cite{Tabuada07}: It turns out that each algebraic
Calabi-Yau category of dimension $2$ containing a cluster-tilting
subcategory is equivalent to a stable derived category
of a differential graded category whose perfect derived
category is Calabi-Yau of dimension~$3$.

A more direct approach, based on the universal
property of the cluster category \cite{Keller05}, has been discovered
by Michel Van den Bergh, who has kindly accepted to include his proof
as an appendix to this article.

It turns out that the main theorem and its proofs can be generalized
almost without effort to Calabi-Yau categories of any dimension $d\geq
2$. However, one has to take into account that in the $d$-cluster
category, the selfextensions of the canonical cluster tilting object
vanish in degrees $-d+2, \ldots, -1$. This condition therefore has to
be added to the hypotheses of the generalized main theorem.

\subsection*{Acknowledgments} This research started during a
stay of the first-named author at the Norwegian University of Science
and Technology (NTNU). He thanks the second-named author
and the members of her group at the NTNU for their warm hospitality.
Both authors thank Michel Van den Bergh for pointing out gaps
and detours in the original proof and for agreeing to include his own
proof as an appendix to this article.
They are grateful to Carl Fredrik Berg for pointing
out reference \cite{Benson98}.

\section{The main theorem and two applications}

\subsection{Statement} \label{subsection:statement}
Let $k$ be a perfect field. Let $\ce$ be a $k$-linear Frobenius
category with split idempotents. Suppose that its stable category
$\cc=\ul{\ce}$ has finite-dimensional $\Hom$-spaces and is
Calabi-Yau of CY-dimension $2$, \ie we have bifunctorial isomorphisms
\[
D\cc(X,Y) \iso \cc(Y,S^2X) \ko X,Y\in\cc\ko
\]
where $D$ is the duality functor $\Hom_k(?,k)$ and $S$ the
suspension of $\cc$.

Let $\ct\subset\cc$ be a cluster tilting subcategory. Recall
from \cite{KellerReiten07} that this means that $\ct$ is a
$k$-linear subcategory which is functorially finite in $\cc$
and such that an object $X$ of $\cc$ belongs to $\ct$ iff
we have $\Ext^1(T,X)=0$ for all objects $T$ of $\ct$.
As shown in \cite{KellerReiten07}, the category $\mod\ct$ of
finitely presented $\ct$-modules is then abelian. If it is hereditary,
the cluster category $\cc_\ct$, as defined in
\cite{BuanMarshReinekeReitenTodorov04},
is the orbit category of the bounded derived category $\cd^b(\mod \ct)$
under the action of the autoequivalence $S^2 \circ \Sigma^{-1}$
where $S$ is the suspension and $\Sigma$ the Serre functor
of $\cd^b(\mod \ct)$.

\begin{theorem}
If $\mod\ct$ is hereditary,
then $\cc$ is triangle equivalent to the cluster category
$\cc_\ct$.
\end{theorem}

We will prove the theorem in section~\ref{section:proofs} below.
Now assume that $k$ is algebraically closed.
Let $\cR$ be the radical of $\ct$, \ie the ideal such that for two
indecomposables $X,Y$, the space $\cR(X,Y)$ is formed by the non
isomorphisms from $X$ to $Y$. Let $Q$ be the quiver of $\ct$: Its
vertices are the isomorphism classes of indecomposables of $\ct$
and the number of arrows from the class of an indecomposable $X$
to an indecomposable $Y$ is the dimension of the vector space
$\cR(X,Y)/\cR^2(X,Y)$.

\begin{corollary}
If $k$ is algebraically closed and
for each vertex $x$ of $Q$, only finitely many paths start
in $x$ and only finitely many paths end in $x$,
then $\cc$ is triangle equivalent to the cluster category
$\cc_Q$.
\end{corollary}

Note that under the assumptions of the corollary, the projective (right)
$kQ$-module $kQ(?,x)$ and the injective $kQ$-module $D kQ(x,?)$
are of finite total dimension and that the category $\mod kQ$
of finitely presented $kQ$-modules coincides with the category
of modules of finite total dimension. We will prove the
corollary in section~\ref{section:proofs} below.

\subsection{Application: Cohen-Macaulay modules}
\label{subsection:CM-modules}
Suppose that $k$
is algebraically closed of characteristic $0$. Let the cyclic group $G=\Z/3\Z$
act on the power series ring $S=k[[X,Y,Z]]$ such that a generator
of $G$ multiplies each indeterminate by the same primitive
third root of unity. Then the fixed point ring $R=S^G$ is a Gorenstein ring,
\cf \eg \cite{Watanabe74}, and an
isolated singularity of dimension $3$, \cf \eg Corollary~8.2 of
\cite{IyamaYoshino06}.
The category $\cmac(R)$ of maximal Cohen-Macaulay modules is
an exact Frobenius category. By Auslander's results
\cite{Auslander76}, \cf Lemma 3.10 of \cite{Yoshino90}, its stable
category $\cc=\ul{\cmac}(R)$
is $2$-Calabi Yau. By work of Iyama \cite{Iyama05}, the module
$T=S$ is a cluster-tilting object in $\cc$. The endomorphism
ring of $T$ over $R$ is the skew group ring $S*G$. Under the
action of $G$, the module $T$ decomposes into three indecomposable direct
factors $T=T_1\oplus T_2\oplus S^G$ and we see that its
endomorphism ring $S*G$ is isomorphic to
the completed path algebra of the quiver
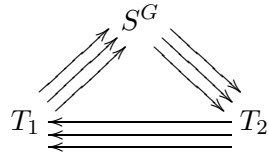
\[
\xymatrix{  & S^G \ar[rd] \ar@<1ex>[rd] \ar@<-1ex>[rd] &  \\
          T_1 \ar[ru] \ar@<1ex>[ru] \ar@<-1ex>[ru]  &
        & T_2 \ar[ll] \ar@<1ex>[ll] \ar@<2ex>[ll] }
\]
subject to all the `commutativity relations' obtained by labelling
the three arrows between any consecutive vertices by $X$, $Y$ and $Z$.
The stable endomorphism ring of $T$ is thus isomorphic to the
path algebra of the generalized Kronecker quiver
\[
Q : \xymatrix{ 1 \ar[r] \ar@<1ex>[r] \ar@<-1ex>[r] & 2}.
\]
The theorem now shows that the stable category of Cohen-Macaulay
modules $\ul{\cmac}(R)$ is triangle equivalent
to the cluster category $\cc_Q$.

    As a further application, we give an alternative proof of a
theorem from \cite{IyamaYoshino06}, stating that the indecomposable
nonprojective rigid modules in $\cmac(R)$ are exactly the
modules of the form $\Omega^i(T_1)$ and $\Omega^i(T_2)$ for $i\in\Z$.
For this, note that the indecomposable rigid objects in $\cc_Q$
are exactly the images of the indeomposable rigid $kQ$-modules and
the $SP$, for $P$ indecomposable projective $kQ$-module
\cite{BuanMarshReinekeReitenTodorov04}. So they
correspond to the vertices of the
component of the AR-quiver of $\cc_Q$ containing the
indecomposable projective $kQ$-modules. The corresponding component of
the AR-quiver of $\ul{\cmac}(R)$
is the one containing $T_1$ and $T_2$.
Hence the indecomposable rigid objects in $\ul{\cmac}(R)$ are all
$\tau$-shifts of these.  Finally, we use that $\tau=\Omega^{-1}$ in this case
\cite{Auslander76}.

Note that this application does not need the full force of the
main theorem. For we only use that the AR-quivers of $\cc_Q$ and
$\ul{\cmac}(R)$ are isomorphic, with the component of the projective
$kQ$-modules for $\cc_Q$ corresponding to the component of $T_1$ and $T_2$
for $\ul{\cmac}(R)$, and it is easy to see that this follows from
proposition 2.1 c) and lemma 3.5 of \cite{KellerReiten07},
\cf also \cite{BuanMarshReiten04}.

\subsection{Application: Non acyclicity} Let $k$ be an algebraically closed
field and $\Lambda$ the preprojective algebra of a simply laced Dynkin
diagram $\Delta$.  Then $\Lambda$ is a finite-dimensional
selfinjective algebra and the stable category $\cc$ of
finite-dimensional $\Lambda$-modules is $2$-Calabi-Yau,
\cf \cite{CrawleyBoevey00}, and admits a
canonical cluster-tilting subcategory $\ct'$ with finitely many
indecomposables, \cf \cite{GeissLeclercSchroeer05b}. Let $Q'$ be its
quiver.  For example, by [loc. cit.], the quivers $Q'$ corresponding
to $\Delta=A_5$ and  $\Delta=D_4$ are respectively
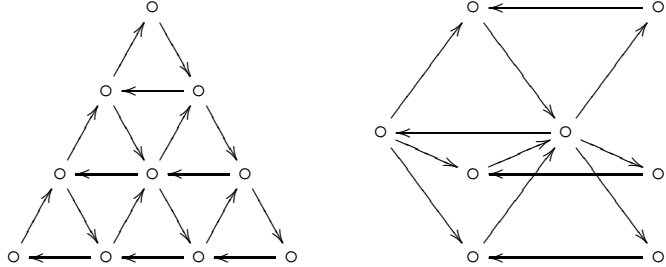
\[
\begin{xy} 0;<0.7pt,0pt>:<0pt,-0.63pt>::
(75,0) *+{\circ} ="0",
(50,50) *+{\circ} ="1",
(100,50) *+{\circ} ="2",
(25,100) *+{\circ} ="3",
(75,100) *+{\circ} ="4",
(125,100) *+{\circ} ="5",
(0,150) *+{\circ} ="6",
(50,150) *+{\circ} ="7",
(100,150) *+{\circ} ="8",
(150,150) *+{\circ} ="9",
"1", {\ar"0"},
"0", {\ar"2"},
"2", {\ar"1"},
"3", {\ar"1"},
"1", {\ar"4"},
"4", {\ar"2"},
"2", {\ar"5"},
"4", {\ar"3"},
"6", {\ar"3"},
"3", {\ar"7"},
"5", {\ar"4"},
"7", {\ar"4"},
"4", {\ar"8"},
"8", {\ar"5"},
"5", {\ar"9"},
"7", {\ar"6"},
"8", {\ar"7"},
"9", {\ar"8"},
\end{xy}
\quad\quad
\begin{xy} 0;<0.7pt,0pt>:<0pt,-0.63pt>::
(0,75) *+{\circ} ="0",
(50,0) *+{\circ} ="1",
(50,100) *+{\circ} ="2",
(50,150) *+{\circ} ="3",
(100,75) *+{\circ} ="4",
(150,0) *+{\circ} ="5",
(150,100) *+{\circ} ="6",
(150,150) *+{\circ} ="7",
"0", {\ar"1"},
"0", {\ar"2"},
"0", {\ar"3"},
"4", {\ar"0"},
"1", {\ar"4"},
"5", {\ar"1"},
"2", {\ar"4"},
"6", {\ar"2"},
"3", {\ar"4"},
"7", {\ar"3"},
"4", {\ar"5"},
"4", {\ar"6"},
"4", {\ar"7"},
\end{xy}
\]
\bigskip
Part b) of the following proposition was obtained independently
by C.~Geiss \cite{Geiss06}.
\begin{proposition} Suppose that $\Lambda$ is representation-infinite.
\begin{itemize}
\item[a)] The stable category $\cc=\ul{\mod} \Lambda$ is not equivalent to the cluster
category $\cc_Q$ of a finite quiver $Q$ without oriented cycles.
\item[b)] The quiver $Q'$ of the canonical cluster-tilting subcategory
of \cite{GeissLeclercSchroeer05b} is not mutation-equivalent to a
quiver $Q$ without oriented cycles.
\end{itemize}
\end{proposition}

In particular, it follows that the two above quivers are not
mutation-equivalent to quivers without oriented cycles. In the proof
of the proposition, we use the main theorem. Let us stress that,
as at the end of \ref{subsection:CM-modules}, we do
not need its full force but only use the isomorphism of AR-quivers.
This is the variant of the proof also given by C.~Geiss \cite{Geiss06}.

\begin{proof} a) Recall first that the AR-translation $\tau$ is isomorphic
to the suspension in any $2$-Calabi-Yau category, so that it is
preserved under triangle equivalences.  We know from
\cite{AuslanderReiten96} that the AR-translation $\tau$ of $\cc$ is
periodic of period dividing $6$.  In particular, we have
$\tau^6(X)\iso X$ for each indecomposable $X$ of $\cc$. But in
$\cc_Q$, for each indecomposable $X$ which is the image of
a preprojective $kQ$-module, the iterated
translates $\tau^{-n}(X)$, $n\geq 0$, are all pairwise non isomorphic
since $Q$ is representation-infinite.

b) Suppose that $Q'$ is mutation-equivalent to a quiver $Q$.
By one of the main results of \cite{GeissLeclercSchroeer07}, it follows that
$\cc$ contains a cluster-tilting subcategory $\ct$ whose quiver is $Q$.
If $Q$ does not have oriented cycles, it follows from the main theorem
that $\cc$ is triangle equivalent to $\cc_Q$ in contradiction to a).
\end{proof}

\section{Proofs} \label{section:proofs}

\subsection{Proof of the corollary}
First recall from \cite{KellerReiten07} that the
category $\mod\ct$ of finitely presented $\ct$-modules is abelian
and Gorenstein of dimension at most $1$. It follows from our
hypothesis that each object of $\mod\ct$ has a finite composition
series all of whose subquotients are simple modules
\[
S_M = \ct(?,M)/\cR(?,M)
\]
associated with indecomposables $M$ of $\ct$ and that each of these
simple modules is of finite projective dimension. Thus each object of
$\mod\ct$ is of finite projective dimension so that $\mod\ct$ has to
be hereditary. Since $k$ is algebraically closed, it follows that
$\mod\ct$ is equivalent to $\mod kQ$ and the claim of the corollary
follows from the theorem.

\subsection{Plan of the proof of the theorem} \label{subsection:plan}
Our aim is to construct a triangle
equivalence $\cc \to\cc_\ct$ such that the triangle
\[
\xymatrix{\cc \ar[r] \ar[rd] & \cc_\ct \ar[d]\\
  & \mod \ct}
\]
becomes commutative, where the diagonal functor takes $X$ to
$\cc(?,X)|\ct$. To
construct the triangle equivalence $\cc \to \cc_\ct$, we use the
construction of $\cc_\ct$ given in \cite{Keller05}, namely, the category
$\cc_\ct$ is the stable derived category $\stab(\ct\oplus D\ct[-3])$ of
the differential graded (=dg) category whose objects are the objects
of $\ct$ and whose morphism complexes are given by the graded modules
\[
\ct(x,y)\oplus (D\ct(y,x))[-3]
\]
endowed with the vanishing differential (the construction of the
stable derived category is recalled in section~\ref{subsection:proof} below).
Thus, we have to construct an equivalence
\[
\cc \to \stab(\ct\oplus (D\ct)[-3]).
\]
We proceed in three steps: 1) We construct a dg category $\ca$ and
a triangle functor
\[
\cc \to \stab(\ca).
\]
We show moreover that the subcategory of indecomposables of the
homology $H^*\ca$ is isomorphic to $\ct\oplus (D\ct)[-3]$.

2) Using the fact that $k$ is perfect we
show that the dg category $\ca$ is formal, \ie linked to its homology
by a chain of quasi-isomorphisms. This yields the required triangle functor
\[
\cc \to\stab(\ca) \iso \stab(H^*\ca) \iso \stab(\ct\oplus
(D\ct)[-3])= \cc_\ct.
\]

3) In a final step, we show that the composed functor
$\cc\to\cc_\ct$ is fully faithful and that its image generates
$\cc_\ct$.

\subsection{The proof}\label{subsection:proof}
Let $\cm\subset\ce$ be the preimage of $\ct$ under the projection
functor. In particular, $\cm$ contains the subcategory $\cp$ of
the projective-injective objects in $\cm$. Note that $\ct$ equals
the quotient $\ul{\cm}$ of $\cm$ by the ideal of morphisms
factoring through a projective-injective.
For each object $M$ of $\ct$, choose an $\ce$-acyclic complex
$A_M$ of the form
\[
0 \to M_1 \to M_0 \to P \to M \to 0 \ko
\]
where $P$ is $\ce$-projective and $M_0, M_1$ are in $\cm$,
\cf \cite{KellerReiten07}. Note that if $\Omega M$ denotes
the kernel of $P\to M$, the induced morphism $M_0 \to \Omega M$
is automatically a right $\cm$-approximation of $\Omega M$.
Let $\ca$ be the dg (=differential graded) subcategory of the dg
category $\cc(\ce)_{dg}$ of complexes over $\ce$ whose objects are
these acyclic complexes. Thus, for two objects $A_L$ and $A_M$ of
$\ca$, we have
\[
H^{n} \ca(A_L,A_M) = \Hom_{\ch\ce}(A_L, A_M[n]) \ko
\]
where $\ch\ce$ denotes the homotopy category of complexes
over $\ce$. To compute this space, let $G_1$ be the functor
$\ce\to \Mod\cm$ taking an object $X$ to $\ce(?,X)|\cm$. The
image of $A_M$ under $G_1$ is a projective resolution of
the $\cm$-module $\ul{M}=\ul{\ce}(?,M)$. Thus we have
\[
\Hom_{\ch\ce}(A_L, A_M[n]) \iso (\cd\cm)(\ul{L}, \ul{M}[n]) =
\Ext^n_\cm(\ul{L},\ul{M}) \ko
\]
where $\cd\cm$ denotes the (unbounded) derived category of
$\Mod\cm$. Notice that by the Yoneda lemma, for each object
$N$ of $\cm$, we have a canonical isomorphism
\[
\Hom_\cm(G_1 N, \ul{M}) \iso \ul{M}(N) = \ul{\ce}(N,M)=\cc(N,M).
\]
Using this we see that the vector space $\Ext^n_\cm(\ul{L},\ul{M})$
is the homology in degree $n$ of the complex
\[
0 \to \cc(L,M)\to 0 \to \cc(L_0,M) \to \cc(L_1, M) \to 0.
\]
Clearly it is isomorphic to $\cc(L,M)$ for $n=0$.
Using the triangle
\[
S^{-2} L \to L_1 \to L_0 \to S^{-1} L
\]
and the fact that $\cc(S^{-1}L_0, M)=0$ and $\cc(S^{-1}L,M)=0$,
we see that the homology is isomorphic to
$\cc(S^{-2}L, M)=D\cc(M,L)$ for $n=3$ and vanishes
for all other $n\neq 0$. More precisely, we see that the map $M \mapsto A_M$
extends to an equivalence whose target is the (additive) graded category
$H^*\ca$ and whose source is the graded category $\ct\oplus
(D\ct)[-3]$ whose objects are those of $\ct$ and whose morphisms
are given by
\[
\ct(L,M) \oplus (D\ct(M,L))[-3].
\]
In particular, we have a faithful functor $\ct \to  H^*\ca$ which
yields an equivalence from $\ct$ to $H^0\ca$. We denote by
$\cd^b(\ca)$ the full subcategory of the derived category $\cd\ca$
whose objects are the dg modules $X$ such that the restriction of
the sum of the $H^nX$, $n\in\Z$,
to $\ct$ lies in the category $\mod\ct$ of finitely
presented $\ct$-modules (by Proposition~2.1 a) of \cite{KellerReiten07},
this category is abelian). In particular, each representable $\ca$-module lies
in $\cd^b(\ca)$ (by Proposition~2.1 b) of \cite{KellerReiten07})
and thus the perfect derived category $\per(\ca)$
is contained in $\cd^b(\ca)$. We denote by $\stab(\ca)$ the
triangle quotient $\cd^b(\ca)/\per(\ca)$. Recall from
\cite{KellerVossieck87} \cite{Rickard89b} that we have a triangle
equivalence
\[
\ul{\ce} \iso \cd^b(\ce)/\cd^b(\cp).
\]
Let $G: \ch^b(\ce) \to \cd \ca$ be the functor which takes
a bounded complex $X$ over $\ce$ to the functor
\[
A_M \mapsto \HOM_\ce(A_M, X) \ko
\]
where $\HOM_\ce$ is the complex whose $n$th component is formed
by the morphisms of graded objects, homogeneous of degree $n$,
and the differential is the supercommutator with the differentials
of $A_M$ and $X$.
We will show that $G$ takes $\cd^b(\cp)$ to zero, that it
maps $\ch^b(\ce)$ to $\cd^b(\ca)$ and the subcategory of acyclic
complexes to $\per(\ca)$. Thus it will induce a triangle functor
\[
\cd^b(\ce)/\cd^b(\cp) \to \cd^b(\ca)/\per(\ca)
\]
and we will obtain the required functor as the composition
\[
\ul{\ce} \iso \cd^b(\ce)/\cd^b(\cp) \to
\cd^b(\ca)/\per(\ca)=\stab(\ca).
\]
First recall that if $A$ is an ayclic complex and $I$ a left
bounded complex of injectives, then each morphism from $A$ to $I$
is nullhomotopic. In particular, the complex $\HOM_\ce(A_M,P)$
is nullhomotopic for each $P$ in $\ch^b(\cp)$. Thus $G$ takes
$\ch^b(\cp)$ to zero. Now, we would like to show that $G$ takes values
in $\cd^b(\ca)$ and that the image of each bounded acyclic complex
is in $\per(\ca)$. For this, we need to compute
\[
(GX)(A_L)=\HOM_{\ce}(A_L, X)
\]
for $L$ in $\cm$ and $X$ in $\ch^b(\ce)$.
To show that the restriction of the sum of the homologies of $GX$ lies in
$\mod\ct$, it suffices to show that this holds if $X$ is concentrated in one
degree. Moreover, if we have a conflation
\[
0 \to M_1 \to M_0 \to X \to 0
\]
of $\ce$ with $M_i$ in $\cm$, it induces a short exact sequence
of complexes
\[
0 \to GM_1 \to GM_0 \to GX \to 0.
\]
So we may suppose that $X$ is an object of $\cm$ considered as a
complex concentrated in degree $0$. Then one computes that the space
\[
\Hom_{\ch\ce}(A_L, X[n])
\]
is isomorphic to the homology in degree $n$ of the complex
\[
0 \to \cc(L,X) \to 0 \to \cc(L_0,X) \to \cc(L_1,X) \to 0 \ko
\]
where $\cc(L,X)$ is in degree $0$. For $n=0$, we find that the
homology is $\cc(L,X)$. Using the triangle
\[
S^{-2} L \to L_1 \to L_0 \to S^{-1}L
\]
and the vanishing of $\cc(S^{-1}L,X)$ and $\cc(S^{-1}L_0,X)$ we see
that the homology in degree $n$ is $\cc(L,S^2 X)$ for $n=3$ and
vanishes for all other $n\neq 0$.
This shows that the restriction of the sum of the
homologies of $GX$ to $\ct$ lies in $\mod\ct$ since the restriction of
$\ul{\ce}(?,Y)$ to $\ct$ lies in $\mod\ct$ for each $Y$ in
$\ul{\ce}$.

Now we have to show that $G$ takes acyclic bounded complexes to
perfect dg $\ca$-modules. For this, we first observe that we have a
factorization of $G$ as the composition
\[
\xymatrix{\ch^b(\ce) \ar[r]^{G_1} & \cd\cm \ar[r]^{G_2} & \cd\ca}
\]
where $G_1$ sends $X$ to $\HOM(?,X)|\cm$ and $G_2$ sends $Y$
to the dg module
\[
A_L \mapsto \HOM(G_1 A_L, Y).
\]
Clearly the
functor $G_1$ sends $\ce$-acyclic bounded complexes to bounded
complexes whose homology modules are in $\mod \ul{\cm}$. Since
$\mod\ul{\cm}$ lies in $\per\cm$, it follows that $G_1$ sends
bounded acyclic complexes to objects of $\per_{\ul{\cm}}(\cm)$.
Under the functor $G_2$, the module $\ul{\cm}(?,L)$ is sent
to $A_L$ and $G_2$ restricted to the triangulated subcategory
generated by the $\ul{\cm}(?,L)$ is fully faithful. We claim
that this subcategory equals $\per_{\ul{\cm}}(\cm)$. Indeed,
each object in $\per_{\ul{\cm}}(\cm)$ is an iterated
extension of its homology objects placed in their respective
degrees. So it suffices to show that each object concentrated in
degree $0$ is the cone over a morphism between objects
$\ul{\cm}(?,L)$, $L\in\cm$. But this is clear since $\mod \ul{\cm}$ is
equivalent to $\mod\ct$, which is hereditary. It follows
that $G_2$ induces an equivalence from $\per_{\ul{\cm}}(\cm)$
to $\per\ca$ and thus $G=G_2 G_1$ sends bounded acyclic
complexes to $\per\ca$.
Thus, we obtain the required triangle functor $F: \cc \to
\stab(\ca)$.

In section~\ref{subsection:formality} below, we will show that
$\ca$ is formal. Thus we get an isomorphism
\[
\ct\oplus (D\ct)[-3] \iso \ca
\]
in the homotopy category of small dg categories. This yields an
equivalence
\[
\cc_\ct = \stab(\ct\oplus (D\ct)[-3]) \iso \stab(\ca).
\]
By construction, it takes each object $T$ of $\ct$ to the module
$T^\wedge=\ct(?,T)$ in $\cc_\ct$. Since $\ct$ generates $\cc$ and
the $T^\wedge$, $T\in\ct$, generate $\cc_\ct$, it is enough to show
that $F$ is fully faithful. We thank Michel Van den Bergh for simplifying
our original argument: For each object $X$ of $\cc$, we
have a triangle
\[
T_1 \to T_0 \to X \to ST_1
\]
with $T_0, T_1$ in $\ct$. Thus, to conclude that $F$ induces
a bijection
\[
\cc(T,X) \to \cc_\ct(FT,FX)
\]
for each $T\in\ct$,
it suffices to show that $F$ induces bijections
\[
\cc(T, T'[i]) \iso \cc_\ct(FT, FT'[i])
\]
for $T,T'$ in $\ct$ and $0\leq i\leq 1$. This is
clear for $i=1$ and not hard to see for $i=0$. We conclude that
for each $Y$ of $\cc$, $F$ induces bijections
\[
\cc(T',Y[i]) \iso \cc_\ct(FT', FY[i])
\]
for all $T'$ in $\ct$ and all $Y$ in $\cc$ and $i\in\Z$. By the
above triangle, it follows that $F$ induces bijections
\[
\cc(X,Y) \to \cc_\ct(FX,FY)
\]
for all $X,Y$ in $\cc$.

\subsection{Formality}\label{subsection:formality} For
categories $\ct$ given by `small enough' quivers $Q$, one can use the
argument of Lemma~4.21 of Seidel-Thomas' \cite{SeidelThomas01} to show
that the category $\ct\oplus (D\ct)[-3]$ is intrinsically formal and
thus $\ca$ is formal. We thank Michel Van den Bergh for pointing out
that for general categories $\ct$ with hereditary module categories,
Seidel-Thomas' argument cannot be adapted.  Instead, we show directly
that $\ca$ is formal (we do not know if $\ct\oplus (D\ct)[-3]$ is
intrinsically formal). Of course, it suffices to show that the full
subcategory $\ca'$ whose objects are the $A_M$ with indecomposable $M$
is formal.

Since $k$ is perfect, the category of bimodules over
a semi-simple $k$-category is still semisimple. From this, one
deduces that the category $\ct$ is
equivalent to the tensor category of a bimodule over the
semi-simplification of $\ct$, \cf Proposition~4.2.5 in
\cite{Benson98}. Using this we can construct a lift of
the functor $\ind\ct\to H^0 (\ca')$ to a functor $\ind\ct \to Z^0(\ca')$,
where $\ind\ct$ denotes the full subcategory of $\ct$ formed
by a set of representatives of the isomorphism classes of
the indecomposables.
We define a $\ct$-bimodule by
\[
X(L,M) = \HOM_\ce(A_L, M) \ko L,M\in \cm \ko
\]
where we consider $M$ as a subcomplex of $A_M$.
Note that $X$ is a right ideal in the category $\ca'$,
that it is a $kQ$-subbimodule of $(L,M) \mapsto \ca'(A_L,A_M)$
and that we have $fg=0$ for all homogeneous elements $f$, $g$
of $X$ of degree $>0$. The computation made above in the proof
that $G$ takes $\ch^b(\ce)$ to $\cd^b(\ca)$ shows that $X$ has
homology only in degree $3$ and that we have a bimodule
isomorphism
\[
D\ct(M,L) \iso H^3 X(L,M) \ko L,M\in \cm.
\]
Thus we have an isomorphism
\[
D\ct[-3] \iso X
\]
in the derived category of $\ct$-bimodules. We choose a projective bimodule
resolution $P$ of $D\ct[-3]$ whose non zero components are
concentrated in degrees $1$, $2$ and $3$ (note that this is possible
since the bimodule category is of global dimension $2$). We obtain
a morphism of complexes of bimodules
\[
P \to X
\]
inducing an isomorphism in homology. We compose it with the inclusion
$X \to \ca'$. All products of elements in the image of $P$ vanish since
they all lie in components of degree $>0$ of $X$. Thus
we obtain a morphism of dg categories
\[
\ct \oplus P \to \ca'
\]
inducing an isomorphism in homology.
This clearly shows that $\ca'$ is formal.

\section{A generalization to higher dimensions}

\subsection{Negative extension groups} \label{subsection:negative}
Let $k$ be a field and  $H$ a finite-dimensional hereditary $k$-algebra.
We write $\nu$ for the Serre functor of the bounded derived
category $\cd=\cd^b(\mod H)$ and $S$ for its suspension functor.
Let $d\geq 2$ be an integer. Let $\cc=\cc^{(d)}_H$ be the {\em $d$-cluster category},
\ie the orbit category of $\cd$ under
the action of the automorphism $\nu^{-1} S^d$, and $\pi : \cd \to \cc$ the
canonical projection functor. We know from \cite{Keller05} that $\cc$
is canonically triangulated and $d$-Calabi Yau and that
$\pi$ is a triangle functor. Moreover, the image $\pi(H)$ of $H$
in $\cc$ is a $d$-cluster tilting object, \cf \eg \cite{KellerReiten07}.
The fact that the module $H$ is projective and concentrated in
degree $0$ yields vanishing properties for the negative selfextension groups
of $\pi(H)$ if $d\geq 3$:

\begin{lemma} We have
\[
\Hom(\pi(H), S^{-i}\pi(H))=0
\]
for $1 \leq i \leq d-2$.
\end{lemma}

\begin{proof} Put $T=H$. For $p\in\Z$,  let $\cd_{\leq p}$ and
$\cd_{\geq p}$ be the $(-p)$th suspensions of the canonical left,
respectively right, aisles of $\cd$, \cf \cite{KellerVossieck88}.
We have to show that the groups
\[
\Hom(T,  \nu^{-p} S^{pd-i} T)
\]
vanish for all $p\in\Z$ and all $1\leq i \leq d-2$. Suppose that $p=-q$ for
some $q\geq 0$. Then we have
\[
\Hom(T,  \nu^{-p} S^{pd-i} T)) = \Hom(T, \nu^{q} S^{-qd-i} T)
\]
and the last group vanishes since $T$ lies in $\cd_{\leq 0}$ and
$\nu^{q} S^{-qd-i}T$ lies in $\cd_{\geq q(d-1)+i}$ and we have
$q(d-1)+i>0$. Now suppose that $p\geq 1$. Then we have
\[
\Hom(T,  \nu^{-p} S^{pd-i} T) = \Hom(\nu^{p} T , S^{pd-i} T) = \Hom(\nu^{p-1}(\nu T), S^{pd-i} T)
\]
and this group vanishes since we have $\nu^{p-1}(\nu T) \in \cd_{\geq p-1}$ (because
$\nu T=\nu H$ is in $\mod H$) and $S^{pd-i}T\in \cd_{\leq -(pd-i)}$ and
\[
(pd-i)-(p-1) = p(d-1) -i+1 \geq d-1-i+1\geq d-i \geq 2.
\]
\end{proof}

\subsection{A characterization of higher cluster categories}
Let $d\geq 2$ be an integer, $k$ an algebraically closed
field and $\cc$ a $\Hom$-finite
algebraic $d$-Calabi Yau category containing a $d$-cluster tilting object $T$.

\begin{theorem} Suppose that $\Hom(T, S^{-i}T)=0$ for $1\leq i\leq d-2$.
If $H=\End(T)$ is hereditary, then, with the notations of section~\ref{subsection:negative},
there is a triangle equivalence $\cc \iso \cc^{(d)}_H$ taking $T$ to $\pi(H)$.
\end{theorem}

Notice that by the lemma above, the assumption on the vanishing of
the negative extension groups is necessary. These assumptions imply
that the endomorphism algebra is Gorenstein of dimension $\leq d-1$,
as we show in lemma~\ref{lemma:Gorenstein} below. For $d\geq 3$,
this does not, of course, imply that the endomorphism
algebra is hereditary if its quiver does not have oriented
cycles, but it implies that the global dimension is
at most $d-1$.

We will prove the theorem below in section~\ref{subsection:proof2}.
In \cite{IyamaYoshino06}, Theorem~1.3, the reader will find an example
from the study of rigid Cohen-Macaulay modules which shows
that the vanishing of the negative extension groups does not follow from
the other hypotheses. The following simple example, based on an idea
of M.~Van den Bergh, is similar in spirit:

\subsection{Example}
Let $\tilde{H}$ be the path algebra of a quiver with underlying
graph $A_6$ and alternating orientation.
Put $\tilde{\cd}=\cd^b(\mod\tilde{H})$ and let $\cc$ be the orbit category of
$\tilde{\cd}$ under the automorphism
$F=\tau^4$ (where $\tau=S^{-1}\nu$). Then $\cc$ is $3$-Calabi Yau: Indeed, one
checks that $F^2=\tau^{-1} S^2$ in $\tilde{\cd}$, which clearly yields $\nu=S^3$ in
$\cc$. The following diagram shows a piece of the Auslander-Reiten
quiver of $\tilde{\cd}$ which is a `fundamental domain' for $F$.
To obtain the Auslander-Reiten quiver of $\cc$, we identify the left and right
borders.
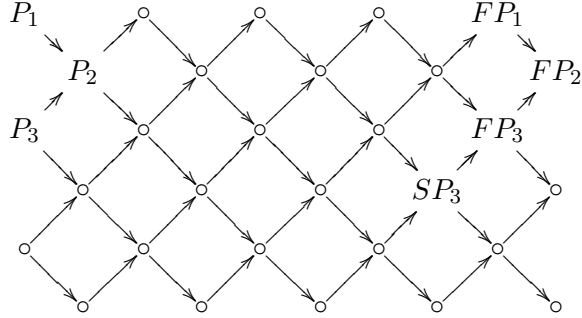
\[
\begin{xy} 0;<1pt,0pt>:<0pt,-1pt>::
(0,0) *+{P_1} ="0",
(22,22) *+{P_2} ="1",
(0,44) *+{P_3} ="2",
(22,67) *{\circ} ="3",
(0,89) *{\circ} ="4",
(22,111) *{\circ} ="5",
(45,0) *{\circ} ="6",
(67,22) *{\circ} ="7",
(45,44) *{\circ} ="8",
(67,67) *{\circ} ="9",
(45,89) *{\circ} ="10",
(67,111) *{\circ} ="11",
(89,0) *{\circ} ="12",
(112,22) *{\circ} ="13",
(89,44) *{\circ} ="14",
(112,67) *{\circ} ="15",
(89,89) *{\circ} ="16",
(112,111) *{\circ} ="17",
(134,0) *{\circ} ="18",
(156,22) *{\circ} ="19",
(134,44) *{\circ} ="20",
(156,67) *+{SP_3} ="21",
(134,89) *{\circ} ="22",
(156,111) *{\circ} ="23",
(179,0) *+{FP_1} ="24",
(201,22) *+{FP_2} ="25",
(179,44) *+{FP_3} ="26",
(201,67) *{\circ} ="27",
(179,89) *{\circ} ="28",
(201,111) *{\circ} ="29",
"0", {\ar"1"},
"2", {\ar"1"},
"1", {\ar"6"},
"1", {\ar"8"},
"2", {\ar"3"},
"4", {\ar"3"},
"3", {\ar"8"},
"3", {\ar"10"},
"4", {\ar"5"},
"5", {\ar"10"},
"6", {\ar"7"},
"8", {\ar"7"},
"7", {\ar"12"},
"7", {\ar"14"},
"8", {\ar"9"},
"10", {\ar"9"},
"9", {\ar"14"},
"9", {\ar"16"},
"10", {\ar"11"},
"11", {\ar"16"},
"12", {\ar"13"},
"14", {\ar"13"},
"13", {\ar"18"},
"13", {\ar"20"},
"14", {\ar"15"},
"16", {\ar"15"},
"15", {\ar"20"},
"15", {\ar"22"},
"16", {\ar"17"},
"17", {\ar"22"},
"18", {\ar"19"},
"20", {\ar"19"},
"19", {\ar"24"},
"19", {\ar"26"},
"20", {\ar"21"},
"22", {\ar"21"},
"21", {\ar"26"},
"21", {\ar"28"},
"22", {\ar"23"},
"23", {\ar"28"},
"24", {\ar"25"},
"26", {\ar"25"},
"26", {\ar"27"},
"28", {\ar"27"},
"28", {\ar"29"},
\end{xy}
\]
Using the mesh category of this quiver,
it is not hard to check that the sum of the images of the indecomposable
projectives $P_1, P_2, P_3$ in $\cc$ is a $3$-cluster tilting
object whose endomorphism ring $H$ is the path algebra on the full subquiver
with the corresponding $3$ vertices. On the other hand,
the image of $P_3$ in $\cc$ has a one-dimensional space of $(-1)$-extensions.
Note that $\cc=\tilde{\cd}/F$ is nevertheless an orbit category and admits
the $3$-cluster category
\[
\cc^{(3)}_{\tilde{H}}=\tilde{\cd}/F^2
\]
as a `2-sheeted covering'.

\subsection{Proof}\label{subsection:proof2}
The proof of the theorem follows the lines of the one in section~\ref{subsection:proof}:
Let $\ct$ be the full subcategory of $\cc$ whose objects are the
direct sums of direct factors of $T$. Let $M$ be an object of $\ct$.
We construct an $\ce$-acyclic complex $A_M$
\[
0 \to M_{d+1} \to M_d \to \ldots \to M_1 \to M_0 \to 0
\]
which yields a resolution of the $\cm$-module
\[
\ul{\ce}(?,M): \cm^{op}\to \Mod k
\]
as in part b) of Theorem~5.4 in \cite{KellerReiten07}. Thus
we can take $M_0=M$ and the morphism $M_1 \to M_0$ is a deflation with projective
$M_1$. Each morphism
\[
M_i \to Z_{i-1}=\ker(M_{i-1} \to M_{i-2}) \ko i\geq 2\ko
\]
yields a $\ct$-approximation in $\cc$. Our vanishing assumption
then implies that $M_2$, \ldots , $M_{d-1}$ are projective.
As in section~\ref{subsection:proof}, we let
$\ca$ be the dg subcategory of the dg
category $\cc(\ce)_{dg}$ of complexes over $\ce$ whose objects are
these acyclic complexes. Thus, for two objects $A_L$ and $A_M$ of
$\ca$, we have
\[
H^{n} \ca(A_L,A_M) = \Hom_{\ch\ce}(A_L, A_M[n]).
\]
One computes that this vector space is isomorphic to $\cc(L,M)$
for $n=0$, to $\cc(L,\Sigma M)=D\cc(M,L)$ for $n=d+1$ and vanishes
for all other $n$. Here we use again our vanishing hypothesis.
We see that the map $M \mapsto A_M$
extends to an equivalence whose target is the (additive) graded category
$H^*\ca$ and whose source is the graded category $\ct\oplus
(D\ct)[-(d+1)]$ whose objects are those of $\ct$ and whose morphisms
are given by
\[
\ct(L,M) \oplus (D\ct(M,L))[-(d+1)].
\]
Now the proof proceeds as in~\ref{subsection:proof} and we obtain
a triangle functor $F: \cc \to \cc_\ct$ taking the subcategory
$\ct$ to $\add \pi(H)$ and whose restriction to $\ct$ is an equivalence.
By lemma~\ref{lemma:equivalences} below, $F$ is an equivalence.

\subsection{Equivalences between $d$-Calabi Yau categories}
Let $d\geq 2$ be an integer, $k$ a field and $\cc$ and $\cc'$
$k$-linear triangulated categories which are $d$-Calabi Yau.
Let $\ct\subset \cc$ and $\ct'\subset \cc'$
be $d$-cluster tilting subcategories. Suppose that $F:\cc\to\cc'$
is a triangle functor taking $\ct$ to $\ct'$.

\begin{lemma} \label{lemma:equivalences}
$F$ is an equivalence iff the restriction of $F$ to $\ct$ is
an equivalence.
\end{lemma}

\begin{proof} It follows from Proposition~5.5 a) of \cite{KellerReiten07},
\cf also part (1) of Theorem~3.1 of \cite{IyamaYoshino06},
that $\cc$ equals its subcategory
\[
\ct*S\ct* \ldots * S^{d-1}\ct
\]
and similarly for $\cc'$. Suppose that the restriction of $F$ to $\ct$
is an equivalence. Let $T\in\ct$. By induction, we see that for each $1\leq i\leq d-1$,
the map
\[
\cc(T,Y) \to \cc'(FT, FY)
\]
is bijective for each $Y\in \ct*S\ct *\ldots * S^i\ct$. Thus the map
\[
\cc(S^j T, Y) \to \cc'(S^j FT, FY)
\]
is bijective for all $j\in\Z$ and $Y$ in $\cc$. Then
it follows that the map
\[
\cc(X,Y) \to \cc'(FX,FY)
\]
is bijective for all $X,Y$ in $\cc$. Thus $F$ is fully faithful.
Since $\ct'$ generates $\cc'$, the functor $F$ is an equivalence.
Conversely, if $F$ is an equivalence and takes $\ct$ to $\ct'$,
then the image of $\ct$ has to be $\ct'$ since $F\ct$ is
maximal $(d-1)$-orthogonal in $\cc'$.
\end{proof}

\subsection{The Gorenstein property for certain $d$-Calabi Yau categories}
\label{subsection:Gorenstein}
Let $d\geq 2$ be an integer, $k$ a field and $\cc$
a $k$-linear triangulated category which is $d$-Calabi Yau.
Let $\ct\subset \cc$ be a $d$-cluster tilting subcategory such that
we have
\[
\cc(T, S^{-i}T')=0
\]
for all $1\leq i \leq d-2$ and all $T$, $T'$ of $\ct$.

\begin{lemma} \label{lemma:Gorenstein}
The category $\mod\ct$ is Gorenstein of dimension less than or
equal to $d-1$.
\end{lemma}

\begin{proof} As in \cite{KellerReiten07}, one sees that the functor
$\Hom(T,?)$ induces an equivalence from the category $S^{d}\ct$ to
the category of injectives of $\mod\ct$. So we have to show that
the $\ct$-module $\Hom(?,S^{d}T)$ is of projective dimension $\leq d-1$
for each $T$ in $\ct$. Put $Y=S^{d}T$. We proceed as in section~5.5 of
\cite{KellerReiten07}: Let $T_0 \to Y$ be a right
$\ct$-approximation of $Y$. We define
an object $Z_0$ by the triangle
\[
Z_0 \to T_0 \to Y \to SZ_0.
\]
Now we choose a right $\ct$-approximation $T_1 \to Z_0$ and define
$Z_1$ by the triangle
\[
Z_1 \to T_1 \to Z_0 \to SZ_1.
\]
We continue inductively constructing triangles
\[
Z_i \to T_i \to Z_{i-1} \to SZ_i
\]
for $1<i\leq d-2$. By proposition~5.5 of \cite{KellerReiten07}, the
object $Z_{d-2}$ belongs to $\ct$.
We obtain a complex
\[
0 \to Z_{d-2} \to T_{d-2} \to \ldots \to T_1 \to T_0 \to Y \to 0.
\]
We claim that its image under the functor
$F:\cc\to \mod\ct$ taking an object $X$ of $\cc$ to
$\cc(?,X)|\ct$ is a projective resolution of $FY$.
Indeed, by induction one checks that the object
$Z_i$ belongs to
\[
S^{d-i-1}\ct * S^{-i}\ct*S^{-i+1}\ct * \ldots * S^{-1}\ct * \ct.
\]
Thus, for each $T\in\ct$, we have $\cc(T, S^{-1} Z_i)=0$ by our vanishing assumptions.
Moreover, the maps $FT_{i+1}\to FZ_i$ are surjective, by construction.
Therefore, the triangle
\[
S^{-1} Z_{i-1} \to Z_i \to T_i \to Z_{i-1}
\]
induces a short exact sequence
\[
0 \to FZ_{i} \to FT_i \to FZ_i \to 0
\]
for each $0\leq i \leq d-2$, where $Z_{-1}=Y$. This implies
the assertion.

\end{proof}

\appendix
\section{An alternative proof of the main theorem}
\begin{center}
\textsc{Michel Van den Bergh}
\end{center}
In this appendix we give a proof of Theorem \ref{subsection:statement}
which is based on the universal property of orbit categories
\cite{Keller05}. We use the same notations as in the main text, but
for the purposes of exposition we will assume that $\Tscr$ consists of
a single object $T$ such that $B=\Cscr(T,T)=kQ$ where $Q$ is a
(necessarily finite) quiver. The extension to more general
$\Tscr$ is routine.

\subsection{The dualizing module}
For use below we recall a version of the Gorensteinness result from
\cite{KellerReiten07}.  Assume that $\Cscr$ is a two-dimensional
$\Ext$-finite Krull-Schmidt Calabi-Yau category with a cluster tilting
object $T$.  Let $B=\Cscr(T,T)$.

For a finitely generated projective right $B$-module we define
$P\otimes_B T$ in the obvious way. For any $M\in \Cscr$, there is a
distinguished triangle (\eg\ \cite{KellerReiten07})
\[
P''\otimes_B T\xrightarrow{\phi} Q''\otimes_B T\r M\r
\]
Now we apply this with $M=T[2]$. Consider a distinguished triangle
\[
 P''\otimes_B T\r Q''\otimes_B T\r T[2]\r
\]
Applying the long exact sequence for
$\Hom_\Cscr(T,-)$ we obtain a corresponding
projective resolution as right module of
the dualizing module of $B$:
\begin{equation}
\label{nice}
0\r P''\r Q''\r DB\r 0
\end{equation}
If we choose any other right module resolution of $DB$
\begin{equation}
\label{nice1}
0\r P'\r Q'\xrightarrow{\alpha} DB\r 0
\end{equation}
then it is equal to \eqref{nice} up to  contractible summands. Hence we obtain
a distinguished triangle
\begin{equation}
\label{nice2}
 P'\otimes_B T\r Q'\otimes_B T\xrightarrow{\alpha'} T[2]\r
\end{equation}
Changing, if necessary, $\alpha'$ by a unit in $B=\End(T[2])$ we
may and we will assume that
$\Hom(T,\alpha')=\alpha$ (under the canonical identifications
$\Hom_\Cscr(T,Q'\otimes_BT)=Q'$ and $\Hom_\Cscr(T,T[2])=DB$).
\subsection{The proof}
We now let $\Cscr$ be as in the main text. By \cite[Thm.\
4.3]{Keller94} we may assume that $\Cscr$ is a strict ($=$ closed
under isomorphism) triangulated subcategory of a derived category
$\cd(\Ascr)$ for some DG-category $\Ascr$. We denote by ${}_B\Cscr$ the
full subcategory of $\cd(B\otimes \Ascr)$ whose objects
are
differential
graded $B\otimes \Ascr$-modules which are in $\Cscr$ when considered as
$\Ascr$-modules. Clearly ${}_B\Cscr$ is triangulated.
\begin{lemmas}
\label{ABlemma}
Assume that $B=kQ$. Then the following holds.
\begin{itemize}
\item[(a)] $T$ may be lifted to an object in ${}_B\Cscr$, also denoted by $T$.
\item[(b)] There is an isomorphism in ${}_B\Cscr$:
$
DB\Lotimes_B T\cong T[2]
$.
\end{itemize}
\end{lemmas}
\begin{proof}
  We may assume that $T$ is a homotopy projective $\Ascr$-module containing a
  summand for each of the vertices of $Q$. Then we may lift the
  action of the arrows in $Q$ to an action of $kQ$ on $T$. Hence (a) holds.

  To prove (b), we choose a resolution of $B$-bimodules
\[
0\r P' \r Q' \xrightarrow{\alpha} DB\r 0
\]
where $P'$, $Q'$ are projective on the right.
Such a resolution may be
obtained by suitably truncating a projective bimodule resolution of $DB$.
 Derived tensoring this
resolution on the right by $T$ and comparing with \eqref{nice2} we find an
isomorphism in $\Cscr$
\begin{equation}
\label{isomorphism}
c:DB\Lotimes_B T\cong T[2]
\end{equation}
between objects in ${}_B\Cscr$. Note that $c$ satisfies $c\circ
\alpha=\alpha'$.

We claim that $c$ in \eqref{isomorphism} is
compatible with the left $B$-actions in $\Cscr$ on both sides.
Let $b\in B$. Then we have a commutative diagram of right $B$-modules
\begin{equation}
\label{tocomparewith}
\begin{CD}
0 @>>> P' @>>> Q' @>>> DB @>>> 0\\
@. @Vb\cdot VV @V b\cdot VV @Vb \cdot VV @.
\\
0@>>> P' @>>> Q' @>>> DB @>>> 0
\end{CD}
\end{equation}
Tensoring on the right by $T$ we obtain a morphism of triangles in $\Cscr$
\[
\begin{CD}
P'\otimes_B T @>>> Q'\otimes_B T @>\alpha'>> T[2] @>>> \\
@V b\cdot VV @V b\cdot VV @V b' VV\\
P'\otimes_B T @>>> Q'\otimes_B T @>\alpha'>> T[2] @>>>
\end{CD}
\]
where $b'=c\circ (b\otimes \Id_T)\circ c^{-1}$.  We need to prove that $b'=b$ under the identification $B=\End_\Cscr(T[2])$.  This
follows easily by applying the functor $\Hom_\Cscr(T,-)$ and comparing to
\eqref{tocomparewith} (using the fact that $\Hom(T,\alpha')=\alpha$).

The proof of (b) can now be completed by invoking the following lemma.
\end{proof}
\begin{lemmas} \label{ll} Assume that $B$ has Hochschild dimension one. Let $M,N\in {}_B\Cscr$. Then the
map
\[
\Hom_{{}_B\Cscr}(M,N)\r \Hom_{\Cscr}(M,N)^B
\]
is surjective (where $(-)^B$ denotes the $B$-centralizer).
\end{lemmas}
\begin{proof}
  Replacing $M$ by a homotopy projective and $N$ by a homotopy injective $B\otimes
  \Ascr$-module one easily obtains the following identity
\[
\RHom_{B\otimes\Ascr}(M,N)=\RHom_{B^e}(B, \RHom_{\Ascr}(M,N))
\]
which yields a spectral sequence
\[
E^{pq}_2: \Ext^p_{B^e}(B, \Ext^q_{\Cscr}(M,N))\Rightarrow \Ext^n_{{}_B\Cscr}(M,N)
\]
Using the fact that $B$ has projective dimension one as bimodule this
yields a short exact sequence
\[
0\r \Ext^1_{B^e}(B, \Ext^{-1}_\Cscr(M,N))\r
\Hom_{{}_B\Cscr}(M,N)\r \Hom_{B^e}(B, \Hom_\Cscr(M,N))\r 0
\]
which gives in particular what we wanted to show.
\end{proof}
\def\Per{\operatorname{Per}}

We can now finish the proof of the main theorem.  By Lemma
\ref{ABlemma}(a) we have a functor
\[
F=-\Lotimes_B T:D^b(B)\r \Cscr
\]
which by \ref{ABlemma}(b) satisfies
\[
F\circ \Sigma[-2]=F\circ (-\otimes_B DB[-2])\cong F
\]
By the universal property of orbit categories \cite{Keller05}, we
obtain an exact  functor
\[
\cd^b(B)/\Sigma[-2]\r \Cscr
\]
which sends $B$ to $T$. This functor is then an equivalence by Lemma
\ref{lemma:equivalences} which finishes the proof.


\begin{thebibliography}{10}

\bibitem{Amiot07}
Claire Amiot, \emph{On the structure of triangulated categories with finitely
  many indecomposables}, arXiv:math/0612141v2 [math.CT], to appear in Bull.
  SMF.

\bibitem{AssemBruestleSchifflerTodorov05}
I.~Assem, T.~Br{\"u}stle, R.~Schiffler, and G.~Todorov, \emph{Cluster
  categories and duplicated algebras}, J. Algebra \textbf{305} (2006), no.~1,
  548--561.

\bibitem{AssemBruestleSchiffler06}
Ibrahim Assem, Thomas Br\"ustle, and Ralf Schiffler, \emph{{Cluster-tilted
  algebras as trivial extensions}}, arXiv:math.RT/0601537, to appear in
  J.~London Math.~Soc.

\bibitem{Auslander76}
Maurice Auslander, \emph{Functors and morphisms determined by objects},
  Representation theory of algebras (Proc. Conf., Temple Univ., Philadelphia,
  Pa., 1976), Dekker, New York, 1978, pp.~1--244. Lecture Notes in Pure Appl.
  Math., Vol. 37.

\bibitem{AuslanderReiten96}
Maurice Auslander and Idun Reiten, \emph{{$D$}{T}r-periodic modules and
  functors}, Representation theory of algebras (Cocoyoc, 1994), CMS Conf.
  Proc., vol.~18, Amer. Math. Soc., Providence, RI, 1996, pp.~39--50.

\bibitem{BuanIyamaReitenScott07}
Aslak Bakke~Buan, Osamu Iyama, Idun Reiten, and Jeanne Scott, \emph{{Cluster
  structures for 2-Calabi-Yau categories and unipotent groups}},
  arXiv:math.RT/0701557.

\bibitem{BuanMarshReinekeReitenTodorov04}
Aslak Bakke~Buan, Robert~J. Marsh, Markus Reineke, Idun Reiten, and Gordana
  Todorov, \emph{Tilting theory and cluster combinatorics}, Advances in
  Mathematics \textbf{204 (2)} (2006), 572--618.

\bibitem{BuanMarshReiten04b}
Aslak Bakke~Buan, Robert~J. Marsh, and Idun Reiten, \emph{{Cluster mutation via
  quiver representations}}, arXiv:math.RT/0412077, to appear in Comm. Math.
  Helv.

\bibitem{BuanMarshReiten06}
Aslak Bakke~Buan, Robert~J. Marsh, and Idun Reiten, \emph{Cluster-tilted
  algebras of finite representation type}, J. Algebra \textbf{306} (2006),
  no.~2, 412--431.

\bibitem{BuanMarshReiten04}
\bysame, \emph{Cluster-tilted algebras}, Trans. Amer. Math. Soc., \textbf{359}
  (2007), no.~1, 323--332, electronic.

\bibitem{BuanMarshReitenTodorov05}
Aslak Bakke~Buan, Robert~J. Marsh, Idun Reiten, and Gordana Todorov,
  \emph{{Clusters and seeds in acyclic cluster algebras}},
  arXiv:math.RT/0510359, to appear in Proc.~AMS, with an appendix by Aslak
  Bakke Buan, Philippe Caldero, Bernhard Keller, Robert Marsh, Idun Reiten and
  Gordana Todorov.

\bibitem{BuanReiten05b}
Aslak Bakke~Buan and Idun Reiten, \emph{{From tilted to cluster-tilted algebras
  of Dynkin type}}, arXiv:math.RT/0510445.

\bibitem{BuanReiten06}
Aslak Bakke~Buan and Idun Reiten, \emph{Acyclic quivers of finite mutation
  type}, Int. Math. Res. Not. (2006), Art. ID 12804, 10.

\bibitem{BuanReitenSeven06}
Aslak Bakke~Buan, Idun Reiten, and Ahmet Seven, \emph{Tame concealed algebras
  and cluster quivers of minimal infinite type}, arXiv:math.RT/0512137, to
  appear in J.~Pure Appl.~Alg.

\bibitem{BaurMarsh06}
Karin Baur and Robert~J. Marsh, \emph{{A geometric description of $m$-cluster
  categories}}, arXiv:math.RT/0607151, to appear in IMRN.

\bibitem{Benson98}
D.~J. Benson, \emph{Representations and cohomology. {I}}, second ed., Cambridge
  Studies in Advanced Mathematics, vol.~30, Cambridge University Press,
  Cambridge, 1998, Basic representation theory of finite groups and associative
  algebras.

\bibitem{BerensteinFominZelevinsky05}
Arkady Berenstein, Sergey Fomin, and Andrei Zelevinsky, \emph{Cluster algebras.
  {III}. {U}pper bounds and double {B}ruhat cells}, Duke Math. J. \textbf{126}
  (2005), no.~1, 1--52.

\bibitem{BurbanIyamaKellerReiten07}
Igor Burban, Osamu Iyama, Bernhard Keller, and Idun Reiten, \emph{Cluster
  tilting for one-dimensional hypersurface singularities}, arXiv:0704.1249v1
  [math.RT].

\bibitem{CalderoChapoton06}
Philippe Caldero and Fr{\'e}d{\'e}ric Chapoton, \emph{Cluster algebras as
  {H}all algebras of quiver representations}, Comment. Math. Helv. \textbf{81}
  (2006), no.~3, 595--616.

\bibitem{CalderoChapotonSchiffler04}
Philippe Caldero, Fr\'ed\'eric Chapoton, and Ralf Schiffler, \emph{Quivers with
  relations arising from clusters (${A}_n$ case)}, Trans. Amer. Math. Soc.
  \textbf{358} (2006), no.~5, 1347--1364.

\bibitem{CalderoKeller05a}
Philippe Caldero and Bernhard Keller, \emph{{From triangulated categories to
  cluster algebras}}, arXiv:math.RT/0506018, to appear in Inv. Math.

\bibitem{CalderoKeller05b}
\bysame, \emph{{From triangulated categories to cluster algebras II}},
  arXiv:math.RT/0510251, to appear in Ann. Scient. ENS.

\bibitem{CalderoReineke06}
Philippe Caldero and Markus Reineke, \emph{On the quiver grassmannian in the
  acyclic case}, arXiv:math/0611074.

\bibitem{CrawleyBoevey00}
William Crawley-Boevey, \emph{On the exceptional fibres of {K}leinian
  singularities}, Amer. J. Math. \textbf{122} (2000), no.~5, 1027--1037.

\bibitem{DerksenWeymanZelevinsky07}
Harm Derksen, Jerzy Weymann, and Andrei Zelevinsky, \emph{Quivers with
  potentials and their representations {I}: {Mutations}}, arXiv:0704.0649v2.

\bibitem{FominZelevinsky02}
Sergey Fomin and Andrei Zelevinsky, \emph{Cluster algebras. {I}.
  {F}oundations}, J. Amer. Math. Soc. \textbf{15} (2002), no.~2, 497--529
  (electronic).

\bibitem{FominZelevinsky03}
\bysame, \emph{Cluster algebras. {II}. {F}inite type classification}, Invent.
  Math. \textbf{154} (2003), no.~1, 63--121.

\bibitem{FominZelevinsky03a}
\bysame, \emph{Cluster algebras: notes for the {CDM}-03 conference}, Current
  developments in mathematics, 2003, Int. Press, Somerville, MA, 2003,
  pp.~1--34.

\bibitem{FominZelevinsky07}
\bysame, \emph{Cluster algebras {IV}: Coefficients}, Compositio Mathematica
  \textbf{143} (2007), 112--164.

\bibitem{Geiss06}
Christof Gei\ss, \emph{Private communication}, 2006.

\bibitem{GeissLeclercSchroeer05b}
Christof Gei{\ss}, Bernard Leclerc, and Jan Schr\"oer, \emph{{Auslander
  algebras and initial seeds for cluster algebras}}, arXiv:math.RT/0506405, to
  appear in Journal of the LMS.

\bibitem{GeissLeclercSchroeer06}
\bysame, \emph{{Partial flag
  varieties and preprojective algebras}}, arXiv:math.RT/0609138.

\bibitem{GeissLeclercSchroeer06a}
\bysame, \emph{Rigid modules over preprojective algebras {II}: The
  {K}ac-{M}oody case}, arXiv:math.RT/0703039.

\bibitem{GeissLeclercSchroeer05c}
\bysame, \emph{{Semicanonical
  bases and preprojective algebras II: A multiplication formula}},
  arXiv:math.RT/0509483, to appear in Compositio Mathematica.

\bibitem{GeissLeclercSchroeer07}
\bysame, \emph{Rigid modules
  over preprojective algebras}, Invent. Math. \textbf{165} (2006), no.~3,
  589--632.

\bibitem{HolmJoergensen06}
Thorsten Holm and Peter Jorgensen, \emph{Cluster categories and
selfinjective
  algebras: type {A}}, arXiv:math/0610728v1 [math.RT].

\bibitem{HolmJoergensen06a}
\bysame, \emph{Cluster categories and selfinjective algebras: type {D}},
  arXiv:math/0612451v1 [math.RT].

\bibitem{Iyama05a}
Osamu Iyama, \emph{Maximal orthogonal subcategories of triangulated categories
  satisfying {S}erre duality}, Oberwolfach Report 6, 2005.

\bibitem{Iyama05}
Osamu Iyama, \emph{Higher dimensional {A}uslander-{R}eiten theory on maximal
  orthogonal subcategories}, Proceedings of the 37th Symposium on Ring Theory
  and Representation Theory, Symp. Ring Theory Represent Theory Organ. Comm.,
  Osaka, 2005, pp.~24--30.

\bibitem{IyamaReiten06}
Osamu Iyama and Idun Reiten, \emph{{Fomin-Zelevinsky mutation and tilting
  modules over Calabi-Yau algebras}}, arXiv:math.RT/0605136, to appear in Amer.
  J. Math.

\bibitem{IyamaYoshino06}
Osamu Iyama and Yuji Yoshino, \emph{{Mutations in triangulated categories and
  rigid Cohen-Macaulay modules}}, arXiv:math.RT/0607736.

\bibitem{Keller94}
Bernhard Keller, \emph{Deriving {D}{G} categories}, Ann. Sci. {\'E}cole Norm.
  Sup. (4) \textbf{27} (1994), no.~1, 63--102.

\bibitem{Keller05}
\bysame, \emph{{On triangulated orbit categories}}, Doc. Math. \textbf{10}
  (2005), 551--581.

\bibitem{KellerReiten07}
Bernhard Keller and Idun Reiten, \emph{{Cluster-tilted algebras are Gorenstein
  and stably Calabi-Yau}}, Advances in Mathematics \textbf{211} (2007),
  123--151.

\bibitem{KellerVossieck87}
Bernhard Keller and Dieter Vossieck, \emph{Sous les cat{\'e}gories
  d{\'e}riv{\'e}es}, C. R. Acad. Sci. Paris S{\'e}r. I Math. \textbf{305}
  (1987), no.~6, 225--228.

\bibitem{KellerVossieck88}
\bysame, \emph{Aisles in derived categories}, Bull. Soc. Math. Belg. S{\'e}r. A
  \textbf{40} (1988), no.~2, 239--253.

\bibitem{MarshReinekeZelevinsky03}
Robert Marsh, Markus Reineke, and Andrei Zelevinsky, \emph{Generalized
  associahedra via quiver representations}, Trans. Amer. Math. Soc.
  \textbf{355} (2003), no.~10, 4171--4186 (electronic).

\bibitem{Rickard89b}
Jeremy Rickard, \emph{Derived categories and stable equivalence}, J. Pure and
  Appl. Algebra \textbf{61} (1989), 303--317.

\bibitem{Ringel90}
Claus~Michael Ringel, \emph{Hall algebras and quantum groups}, Invent. Math.
  \textbf{101} (1990), no.~3, 583--591.

\bibitem{Ringel07}
\bysame, \emph{Some remarks concerning tilting modules and tilted algebras.
  {Origin. Relevance. Future.}}, Handbook of Tilting Theory, LMS Lecture Note
  Series, vol. 332, Cambridge Univ. Press, Cambridge, 2007, pp.~49--104.

\bibitem{SeidelThomas01}
Paul Seidel and Richard Thomas, \emph{Braid group actions on derived categories
  of coherent sheaves}, Duke Math. J. \textbf{108} (2001), no.~1, 37--108.

\bibitem{Tabuada07}
Gon\c{c}alo Tabuada, \emph{{On the structure of Calabi-Yau categories with a
  cluster tilting subcategory}}, Doc. Math. \textbf{12} (2007), 193--213.

\bibitem{Thomas06}
Hugh Thomas, \emph{{Defining an $m$-cluster category}}, arXiv:math.RT/0607173.

\bibitem{Watanabe74}
Keiichi Watanabe, \emph{Certain invariant subrings are {G}orenstein. {I},
  {II}}, Osaka J. Math. \textbf{11} (1974), 1--8; ibid. 11 (1974), 379--388.

\bibitem{Yoshino90}
Yuji Yoshino, \emph{Cohen-{M}acaulay modules over {C}ohen-{M}acaulay rings},
  London Mathematical Society Lecture Note Series, vol. 146, Cambridge
  University Press, Cambridge, 1990.

\bibitem{Zelevinsky04}
Andrei Zelevinsky, \emph{{Cluster algebras: notes for 2004 IMCC (Chonju, Korea,
  August 2004)}}, arXiv:math.RT/0407414.

\bibitem{Zhu06}
Bin Zhu, \emph{{Generalized cluster complexes via quiver representations}},
  arXiv:math.RT/0607155.

\end{thebibliography}

\def\cprime{$'$}
\providecommand{\bysame}{\leavevmode\hbox to3em{\hrulefill}\thinspace}
\providecommand{\MR}{\relax\ifhmode\unskip\space\fi MR }
\providecommand{\MRhref}[2]{%
  \href{http://www.ams.org/mathscinet-getitem?mr=#1}{#2}
}
\providecommand{\href}[2]{#2}

\end{document}